\documentclass{article}
\def \Z {{\mathbf Z}}
\def \R {{\mathbf {R}}}
\def \N {{\mathbf {N}}}
\def \B {{\cal B}}

\usepackage[T2A]{fontenc}
\usepackage[cp1251]{inputenc}
\usepackage[english,russian]{babel}
\usepackage[tbtags]{amsmath}
\usepackage{amsfonts,amssymb}

\begin{document}
%\udk{517.9}
\date{vryzh@mail.ru }

\author{В.\,В.~Рыжиков}
%\address{Московский государственный университет \\
%имени М.В.Ломоносова}
%\email{vryzh@mail.ru}
        
\title{ Множества перемешивания  и жесткость  преобразований}
\markboth{В.\,В.~Рыжиков}{Множества перемешивания  и жесткость}

\maketitle

%\begin{fulltext}

\begin{abstract}{В заметке доказано, что  для всякого бесконечного    
множества $M\subset\N$ нулевой плотности найдется жесткое сохраняющее меру  
 преобразование,  перемешивающее вдоль $M$, что дает   ответ  
на вопрос В. Бергельсона.   
В качестве подходящих примеров  преобразований предлагаются 
    гауссовские действия и пуассоновские надстройки 
над специально подобранными бесконечными конструкциями ранга один.  Эквивалентная формулировка этого результата такая: найдется нормированная мера $\sigma$
на единичной окружности, коэффициенты Фурье которой стремятся к нулю вдоль заданного множества $M$ нулевой плотности
а вдоль некоторого множества вне $M$  стремятся к 1. Эта задача, как мы показываем,  решается методами эргодической теории. 
 Обсуждаются аналоги результата для групповых действий и метод, не использующий гауссовские и  пуассоновские надстройки.

Библиография: 8 названий.}
\end{abstract}
%\begin{keywords}
%    Сохраняющие меру преобразования,  
%мягкое перемешивание, жесткость, перемешивание вдоль множества, действия ранга один, гауссовское действие, пуассоновская %надстройка.
%\end{keywords}

\section{Введение}
Иерархия  перемешивающих свойств сохраняющих меру преобразований включает
$K$-перемешивание, кратное перемешивание, перемешивание, частичное,  легкое (light), мягкое (mild) и слабое   перемешивание (см., например,  \cite{KSF},\cite{A}).
Перемешивание, мягкое перемешивание и  слабое перемешивание  
  отражают  свойства  унитарных операторов, индуцированных  преобразованиями.
Кроме перечисленных  инвариантов имеет место   большое разнообразие асимптотических свойств спектральной природы, в определении которых фигурируют  подмножества времени динамической системы. 
Дадим необходимые для дальнейшего определения.

 Пусть обратимое преобразования $T:X\to X$ 
 вероятностного пространства $(X,\B,\mu)$ сохраняет меру $\mu$. \it 
Перемешивание вдоль множества $M\subset\N$ \rm  означает, что для всех  $A,B\in\B$ выполнено 
 $$\mu ( T^{m}A\cap B)\to\mu(A)\mu(B), \ \  m\to\infty, \ m\in M. $$
 
При  $M=\N$ такое преобразование $T$ называется \it перемешивающим.\rm

Говорят, что сохраняющее меру обратимое преобразование $T$ \it слабо перемешивает, \rm
если оно перемешивает вдоль некоторого множества. Хорошо известно, что  
для всякого слабо перемешивающего преобразования найдется множество плотности 1, вдоль которого преобразование перемешивает.

\it Мягкое перемешивание \rm означает следующее:  для любых   $A, \mu(A)>0,$ и  $h_j\to\infty$
сходимость  $\mu (A\Delta T^{h_j}A)\to 0$ влечет за собой $A=X$.\rm

Если  для некоторой последовательности $h_j\to\infty$ для всех $A\in \B$ выполнено $\mu (A\Delta T^{h_j}A)\to 0$, 
 такое преобразование $T$ называется \it жестким.\rm  Жесткое преобразование, очевидно, не обладает мягким перемешиванием.

Вопрос  В.Бергельсона: 
\it найдется ли   множество $M\subset\N$ нулевой плотности такое,
что  всякое  преобразование $T$,  перемешивающее вдоль $M$,   является мягко перемешивающим?\rm

 Предлагается следующее решение. Для множества $M$ нулевой плотности строим бесконечное преобразование $R$ ранга один, которое,  являясь жестким,   перемешивает\footnote{Для бесконечных преобразований 
свойство перемешивания сохраняет спектральный смысл, но, вообще говоря, теряет физический.}   вдоль множества $M$:
 $$\mu ( R^{m}A\cap B)\to 0, \ \  m\to\infty, \ m\in M, $$ 
для всех  $A,B$ конечной меры.
Рассматривая гауссовскую или пуассоновскую надстройку над $R$  приходим к следующему утверждению:
 \it для любого множества $M\subset \N$ нулевой плотности найдется жесткое сохраняющее меру преобразование
вероятностного пространства, которое  перемешивает вдоль $M$.  \rm

Тем самым мы даем отрицательный ответ на упомянутый вопрос.
Следующее утверждение содержит простой  рецепт приготовления  перемешивающих множеств 
ненулевой плотности для жестких преобразований.

\vspace{3mm}
\bf Теорема 1. \it 
Если  дополнение к множеству $\bar M \subset \N$  для некоторых последовательностей  $a_i,L_i\to \infty$
содержит  объединение    наборов    интервалов  вида 
$$ [a_i, a_i+L_i],  [2a_i, 2a_i+L_i], \dots, [ia_i, ia_i+L_i],$$
то $\bar M$ является перемешивающим множеством для некоторого жесткого преобразования.
  \rm

\vspace{3mm}
Понятно, что 	множество нулевой плотности всегда можно поместить в указанное  множество $\bar M$.
Если  в аналогичной задаче  потребовать  отсутствие  перемешивания, не требуя  жесткость, 
 структура дополнения к множеству $\bar M$ упрощается.

%\newpage
\bf Теорема 2. \it Множество 
$\bar M$ является перемешивающим множеством для некоторого неперемешивающего  преобразования,
если  дополнение к $\bar M$  содержит интервалы $ [a_i, a_i+L_i]$ для некоторых последовательностей  $a_i,L_i\to \infty.$  \rm

\vspace{3mm}
\bf Благодарности. \rm Автор признателен Эль Абдалауи (El Abdalaoui), привлекшего   внимание  к  одному из
 вопросов по  тематике заметки, и особенно   В. Бергельсону (V. Bergelson) и Ж.-П. Тувено (J.-P. Thouvenot) 
за полезные обсуждения.

\section{Инструментарий}

 Позже нам понадобятся   преобразования ранга один,  заданные 
параметрами   
 $$s_j(1)=0, s_j(2)=0,\dots, s_j(r_j-1)=0,  s_j(r_j)=s_j, \ \ s_j\to\infty$$
(похожие примеры преобразований, но  для других целей, использовались   в заметке \cite{19}).
Напомним  смысл этих параметров (наша дальнейшая цель --
 найти подходящую последовательность $s_j$  для заданного множества $\bar M$).     

\bf Конструкции ранга один. \rm 
 Фиксируем натуральное число $h_1$ (высота башни на этапе 1), последовательность $r_j$ 
($r_j\geq 2$ -- число колонн, на которое разрезается  башня на этапе $j$)
и последовательность  целочисленных наборов 
$$ \bar s_j=(s_j(1), s_j(2),\dots, s_j(r_j-1),s_j(r_j)), \ s_j(i)\geq 0.$$ 
Параметры $r_j$, $\bar s_j$ полностью  определяют конструкцию преобразования ранга 1, 
которая строится поэтапно.

Пусть на  шаге $j\geq 1$ определена 
система  непересекающихся полуинтервалов, одинаковой длины 
$$B_j,\ RB_j, \ R^2B_j,\dots, R^{ h_j-1}B_j,$$
причем на  полуинтервалах  $B_j, RB_j, \dots, R^{ h_j-2}B_j$
пребразование $R$ действует как   перенос полуинтервалов.  Такая система называется башней высоты $h_j$.

 Представим 
 $B_j$ в виде  дизъюнктного  объединения  полуинтервалов $ B_j^i$, $1\leq i \leq r_j$, одинаковой длины. Набор  
$$B_j^i, RB_j^i ,R^2 B_j^i,\dots, R^{h_j-1}B_j^i$$ называется 
$i$-ой колонной  этапа $j$.
Над колонной надстраиваем  $s_j(i)$ непересекающихся полуинтервалов меры $\mu(B_j^i)$, 
 получая набор 
$$B_j^i,\ RB_j^i,\ R^2 B_j^i,\dots, R^{h_j-1}B_j^i,\ R^{h_j}B_j^i, \ R^{h_j+1}B_j^i, \dots, R^{h_j+s_j(i)-1}B_j^i.$$

  Терерь соберем все   надстроенные колонны в одну башню. Для этого положим 
$$R^{h_j+s_j(i)}B_j^i = B_j^{i+1}$$ при $i<r_j.$
Обозначая  $B_{j+1}= E^1_j$, получили  башню этапа $j+1$:
$$B_{j+1}, \ RB_{j+1},\  R^2 B_{j+1},\dots, R^{h_{j+1}-1}B_{j+1},$$
где 
$$ h_{j+1} =h_jr_j +\sum_{i=1}^{r_j}s_j(i).$$

Продолжая построение до бесконечности,  получим обратимое преобразование $R:X\to X$,  
 где  $X$ -- объединение всех рассматриваемых полуинтервалов. 
Преобразование $R$   сохраняет меру Лебега. 

\vspace{3mm}
\bf  Вероятностное пространство  Пуассона.  \rm Пусть пространство $X$ как объединение интервалов имеет бесконечную меру. Рассмотрим    пространство конфигураций $X_\ast$, состоящее из всех бесконечных 
  счетных множеств $x_\ast$  таких, что в каждом интервале   пространства $X$ находится лишь конечное число элементов множества $x_\ast$.
Для измеримого подмножества $A\subset X$ конечной $\mu$-меры  определим подмножества конфигураций  
$C_{A,k}$, $k=0,1,2,\dots$,
 в  $X_\ast$ формулой
$$C_{A,k}=\{x_\ast\in X_\ast \ : \ |x_\ast\cap A|=k\}.$$

  Всевозможные конечные пересечения вида $\bigcap_{i=1}^N C(A_i,k_i)$    образуют полукольцо. 
На этом полукольце определим  меру $\mu_\ast$. Потребуем следующее:
если множества $A$, $B$ не пересекаются,  то 
вероятность $\mu_\ast(C_{A,k}\bigcap C_{B,n})$  одновременного пребывания $k$ точек конфигурации  $x_\ast$ в $A$ 
и  $n$ точек  конфигурации  $x_\ast$ в $B$ по определению равна произведению   $\mu_\ast(C_{A,k})\mu_\ast( C_{B,n})$
(события $C_{A,k}$ и $C_{B,n}$ независимы).
 Для набора измеримых непересекающихся  множеств $A_1, A_2,\dots, A_N\subset X$  конечной $\mu$-меры  положим
$$\mu_\ast(\bigcap_{i=1}^N C_{A_i,k_i})=\prod_{i=1}^N \frac {\mu(A_i)^{k_i}}{k_i!}e^{-\mu(A_i)}.$$
 Мера $\mu_\ast$  продолжается с полукольца  на
  вероятностное пространство  $(X_\ast,\B_\ast,\mu_\ast)$.

% \newpage
\bf Пуассоновское и гауссовское действие.  \rm Пусть $R$ -- бесконечное преобразование. Надстройка $R_\ast$
действует на точки в $X_\ast$ (неупорядоченные счетные наборы $\{x_1,x_2,\dots,\}$)  по правилу
 $$R_\ast \{x_1,x_2,\dots,\}= \{Rx_1,Rx_2,\dots,\}.$$ 
Она сохраняет меру $\mu_\ast$ и наследует асимптотические свойства: например, если $\mu(R^{m_i}A\cap B)\to 0$
для всех $A,B$ конечной меры, то $$\mu_\ast ( R_\ast^{m_i}C\cap D)\to\mu(C)\mu(D)$$ для всех  $C,D\in \B_\ast$.
Надстройка  $R_\ast$ также наследует свойство жесткости бесконечного преобразования $R$.

  Преобразованию $R$ соответствует  ортогональный оператор в $L_2(\mu)$, 
которому отвечает гауссовский автоморфизм $G$ пространства $(\R^\infty, \gamma^\infty)$, где  
$\gamma$  -- гауссовская мера на $\R$. Как известно, автоморфизм $G$ спектрально изоморфен 
пуассоновской надстройке $R_\ast$ 
и обладает одинаковыми с ней    асимптотическими свойствами жесткости и перемешивания вдоль множеств. 
Подробнее о  пуассоновских надстройках и гауссовских действиях см., например, \cite{KSF},\cite{N}.

\section{Конструкция бесконечного преобразования}
\bf  Доказательство теоремы 1. \rm Фиксируем множество $M$  нулевой плотности. Положим $r_j=j$.
Будем искать  такую  последовательность $s_j$, чтобы   конструкция $R$ с параметрами   
 $$s_j(1)=0, s_j(2)=0,\dots, s_j(j-1)=0,  s_j(j)=s_j, \ \ r_j\to\infty$$
удовлетворяла условиям теоремы.

Будем предполагать выполненым условие $s_j>>h_j$, поэтому  пространство $X$   имеет бесконечную меру.
Свойство перемешивания вдоль $\bar M$ для бесконечных преобразований означает сходимость
$$\mu ( R^{m}A\cap B)\to 0, \ \ m\to\infty,  \ m\in \bar M,$$
для всех множеств $A,B$ конечной меры.

Предположим, что все параметры на этапах до $j-1$ (включительно) уже заданы, они  определяют значение $h_j$.
Находим такое число $h_{j+1}$,  чтобы при $n=1,2,\dots,j$ все отрезки $[nh_{j+1}-jh_j,nh_{j+1}+jh_j]$
 не пересекались с множеством $\bar M$. Обозначим объединение этих отрезков через $V_j$. 
Положим $s_j= h_{j+1}- jh_j$. 
Последовательность $s_j$ определена  и вместе с ней задана  конструкция $R$.
Преобразование $R$  обладает жесткостью: для всех $A, \,\mu(A)<\infty,$
$$\mu (A\Delta R^{h_j}A)\to 0.$$
Действительно,  если множество $A$ состоит из  этажей некоторой башни,
то для  достаточно больших значений $j$ имеем
$$ \mu (A\cap R^{h_j}A)= \mu(A)-\frac{ \mu(A)} j $$
(по построению $R^{h_j}A$ совпадает с  $A$ на всех колоннах, кроме первой, с которой его пересечение пусто).

Пусть  множества $A,B$ состоят из этажей некоторой башни и для больших $j$ выполнено 
$$ h_j\leq p<h_{j+1}, \ \ \mu(R^p A\cap B)>0.$$
Тогда  $p$ попадает в один из отрезков $[nh_{j+1}-jh_j,nh_{j+1}+jh_j]$, лежащий по построению 
в объединении множеств $V_j.$ 
Следовательно, при $m\in \bar M$ для больших значений $m$ имеем  $\mu(R^m A\cap B)=0,$ 
значит, бесконечная жесткая конструкция $R$ перемешивает вдоль множества $\bar M$.

Чтобы получить нужное преобразование на вероятностном пространстве, 
рассмотрим пуассоновскую  надстройку
$R_\ast$ над $R$. Сохраняя вероятностную меру $\mu_\ast$, надстройка $R_\ast$ 
 наследует свойство жесткости и перемешивание вдоль множества.

Действительно, из $ \mu \left(A\Delta R^{h_j}A\right)\to 0$ вытекает
$$ \mu_\ast \left(C_{A,k}\bigcap R_\ast^{h_j}C_{A,k}\right)=\mu_\ast 
\left(C_{A,k}\bigcap C_{R^{h_j}A,k}\right)\to\mu_\ast(C_{A,k}),$$
а  условие $\mu ( R^{m_j}A\cap B)\to 0$ влечет за собой 
 $$\mu_\ast\left(R_\ast^{m_j}C_{A,k}\bigcap C_{B,n}\right)=\mu_\ast\left(C_{R^{m_j}A,k}\bigcap C_{B,n}\right)
\to \mu_\ast(C_{A,k}) \mu_\ast( C_{B,n}).$$
Требуемое  преобразование построено. Теорема 1 доказана.

\bf Доказательство теоремы 2. \rm Схема доказательства аналогична  схеме доказательства теоремы 1,
но конструкция преобразования $R$ теперь значительно проще: она задается параметрами
$$ r_j=2,  \ s_j(1)=0, \ s_j(2)=s_j.$$

Подбирая быстро растущую последовательность $s_j$, получаем, что вне интервалов $[h_j-h_{j-1}, h_j+h_{j-1}]$,
которые по построению лежат в объединении наперед заданных интервалов $ [a_i, a_i+L_i]$, преобразование $R$ перемешивает.
Также оно  обладает  свойством частичной жесткости: для всех $A$ конечной меры выполнено
$$\mu(R^{h_j}A\cap A)\to \mu(A)/2.\eqno (\ast)$$
Пуассоновская надстройка $R_\ast$  удовлетворяет условиям теоремы 2. 
В \cite{19}  доказательство Утверждения 1 фактически содержит  доказательство  следующей теоремы, по смыслу близкой 
к теореме 2.

\vspace{3mm}
\bf Теорема 3. \it Для любого бесконечного множества  
$\bar H\subset \N$ найдется такое его подмножество $\{h_j: j\in \N\}$, что  $h_j$  являются высотами бесконечной конструкции
$R$ с параметрами $ r_j=2,  \ s_j(1)=0, \ s_j(2)=h_{j+1}-2h_j$,  причем для всяких множеств  $A,B$ конечной меры  условие 
$$\forall k \ \ \mu (R^{n(k)}A\cap B)> 0$$
влечет за собой
$$n(k)=  h_{j_1(k)} \pm  h_{j_2(k)}\dots \pm  h_{j_p(k)} + s,\eqno (sp) $$ 
где $s=s(k)$ и  $p=p(k)$  -- ограниченные последовательности, и  $j_1(k)>j_2(k)>\dots >j_p(k)$.
\rm

\vspace{3mm}
Таким образом,  для подходящей конструкции  $R$ мы формируем 
  основную неперемешивающую последовательность $h_j\to\infty$ такую,
что  все остальные неперемешивающие последовательности располагаются рядом с ней.  Здесь подразумевается то, что
очень быстрый рост последовательности $h_j$ обеспечивает выполнение соотношения
$$  h_{j_1(k)} >>  h_{j_2(k)}+\dots +  h_{j_p(k)}.$$
Пуассоновская надстройка $R_\ast$ (и соответствующий гауссовский автоморфизм)
  наследует перемешивающие свойства преобразования $R$.  Всякая последовательность $m_i\to +\infty$, не имеющая бесконечной подпоследовательности  вида $(sp)$, является перемешивающей.

\section{Заключительные замечания}
Бесконечные конструкции $R$ с параметрами $ r_j=2,  \ s_j(1)=0, \ s_j(2)=s_j$ рассматривались в \cite{19}.
Если $s_j>>h_j$, то для надстройки $R_\ast$ выполнено необычное свойство: \it
 спектр  степеней   $R_\ast^m\otimes R_\ast^m$, $m\geq 1,$ сингулярен в силу свойства $(\ast)$, 
но для всех $n> m\geq 1$ спектр  произведения $R_\ast^m\otimes R_\ast^n$  
содержит  счетно-кратную лебеговскую компоненту.\rm

Конструкции  с параметрами $$s_j(1)=0, s_j(2)=0,\dots,   s_j(j)=s_j, \ \ s_j>>h_j$$
приводит к построению потока $T_t$ на вероятностном пространстве, для которого
для всех рациональных $a$   произведение $ T_t\otimes T_{at}$ является жестким преобразованием, следовательно,
имеет сингулярный спектр, но для почти всех $a$ такое произведение  имеет счетно-кратную лебеговскую компоненту в спектре.

В   \cite{11} изучались  неперемешивающие  конструкции конечной меры, 
для которых осуществлялся  контроль перемешивания вне множеств вида $(sp)$.   
Ниже, сохраняя  контроль  перемешивания, мы  модифицируем эти  конструкции, наделяя их   жесткостью.
Отметим, что    для всех $\kappa \in (0,1)$ конструкции обладает также свойством $\kappa$-перемешивания: 
одновременного сочетания жесткости и перемешивания вдоль последовательности. Напомним, что  для преобразования $T$  вероятностного пространства $\kappa$-перемешивание означает  
сходимость  
$$\forall A,B\in\B  \ \ \mu(T^{k_i}A\cap B)\to \kappa \mu(A)\mu(B) +(1-\kappa)\mu(A\cap B)$$
для некоторой последовательности $k_i\to\infty$.
Это свойство  играет важную роль в спектральной теории динамических систем, так как  оно влечет за собой дизъюнктность сверток спектральной меры преобразования, см. \cite{O69},\cite{S86}.

\vspace{3mm}
\bf Класс жестких конструкций, случай конечной меры. \rm  Пусть $r_j\to\infty$, для большинства (перемешивающих) этапов $j$ параметры  $s_j(i)$  имет вид лестницы:     
$$ s_j(1)=1, \ \ s_j(2)=2,\dots, s_j(j-1)=j-1, \ \  s_j(j)=j. $$ 
Для оставшихся (жестких) этапов $j$ пусть $r_j=N_j q_j$, $q_j\to\infty$, а параметры  $s_j(i)$ повторяют  
лестничную надстройку  $N_j$ раз:
 $$ s_j(1)=1, \ \ s_j(2)=2,\dots \   \dots \    \dots \   \dots \   \dots \   \dots \   \dots \ \dots\   \dots, s_j(q_j)=q_j,$$          
 $$ s_j(q_j+1)=1,\ \ \ s_j(q_j+2)=2,\dots \ \dots \    \dots \    \dots \ \dots\ \dots,  s_j(2q_j)=q_j,$$
$$ \dots  \  \dots\ \dots  \  \ \ \  \dots \    \dots \   \dots \ \ \ \ \dots \   \dots \   \dots $$
$$s_j(N_jq_j-q_j+1)=1, \ s_j(N_jq_j-q_j+2)=2, \dots \   \dots, \ s_j(N_jq_j)=q_j,$$
$$s_j(N_jq_j+1)=1,\ \ \ \  s_j(N_jq_j+2)=2,\dots\dots \dots, s_j((N_j+1)q_j)=q_j.$$ 

%\begin{center}
%\includegraphics{stair}
%\end{center}
В этом классе при  $N_j\to\infty$  можно найти подходящую конструкцию конечной меры, перемешивающую вдоль 
заданного множества $\bar M$, удовлетворяющего условию теоремы 1.   
Чем меньше дополнение к множеству $\bar M$, тем реже жесткие этапы.   Неперемешивающие последовательности асимптотически 
группируются вокруг чисел $n_j(h_j+ 1+2+\dots+q_j)$, где $0\leq n_j <N_j$ и такие  группы можно разместить в наперед заданных
интервалах, дождавшись, когда  длины интервалов станут достаточно большими. 
Контроль перемешивания использует оригинальную технику Адамса, адаптированную к нашей ситуации в духе работы  \cite{11}.

\bf Жесткость   $\Z^n$-действий и множества перемешивания.\rm 
В случае $\Z^n$-действий   несложно реализовать идею построения  подходящей конструкции ранга один (см., например,  в \cite{17}) с   основной неперемешивающей последовательностью, вдоль которой располагаются 
  все остальные неперемешивающие последовательности. 

\vspace{3mm}
\bf Теорема 4.  \it Для любого множества $M\subset \Z^n$ нулевой плотности найдется жесткое действие  группы $\Z^n$,
которое  перемешивает вдоль $M$.  \rm

\vspace{3mm}
Можно пополнить список групп, для которых имеет место похожее утверждение,
однако в общем случае остается актуальной интересная постановка задачи.

%\hfill vryzh@mail.ru

\end{document}